
\baselineskip=14pt
\parskip=10pt
\def\Tilde{\char126\relax}
\def\halmos{\hbox{\vrule height0.15cm width0.01cm\vbox{\hrule height
 0.01cm width0.2cm \vskip0.15cm \hrule height 0.01cm width0.2cm}\vrule
 height0.15cm width 0.01cm}}
\font\eightrm=cmr8  
\font\eighttt=cmtt8
\magnification=\magstephalf

\parindent=0pt
\overfullrule=0in
\bf
\centerline
{REVEREND CHARLES to the aid of MAJOR PERCY and FIELDS-MEDALIST ENRICO}
\rm
\bigskip
\centerline{ {\it Doron ZEILBERGER}\footnote{$^1$}
{\eightrm  \raggedright
Department of Mathematics, Temple University,
Philadelphia, PA 19122, USA. 
{\eighttt E-mail:zeilberg@math.temple.edu; 
WWW: http://www.math.temple.edu/\Tilde zeilberg;
Ftp: ftp.math.temple.edu, directory /pub/zeilberg.}
\break
Supported in part by the NSF. July 21, 1995. To appear in the Amer. Math. 
Monthly.
} 
}
Voltaire said that Archimedes had more imagination than Homer.
Unfortunately, as far as we know, Archimedes's only use of it,
outside of mathematics, was to condemn goldsmiths and to kill people.
More peaceful uses of mathematicians' imagination, to  the outside of
mathematics, in decreasing order of impact, were provided by
Multi-Millionaire Richard Garfield, the Reverend Charles
Dodgson, and Major Percy MacMahon,
who respectively developed : `Magic: The Gathering' ({\it the} game of our
decade), Alice, and an earlier version of Instant Insanity.
 
This is not to say that their imagination did not also help mathematics
proper.
In this {\it quickie}, I observe how Dodgson's[D] rule for
evaluating determinants: (For any $n \times n$ matrix $A$,
let $A_r(k,l)$ be the $r \times r$ submatrix whose upper leftmost
corner is the entry $a_{k,l}$,)
$$
\det A \,\,=\,\,
{{
\det A_{n-1}(1,1) \det A_{n-1}(2,2) 
-
\det A_{n-1}(1,2) \det A_{n-1}(2,1) }
\over
{\det A_{n-2}(2,2)}} \quad ,
$$
immediately implies MacMahon's[M] determinant evaluation:
$$
\det \left [ {{a+i} \choose {b+j}}_{1 \leq i,j \leq n}\right ]
=
{{(a+n)!!(n-1)!!(a-b-1)!!(b)!!} \over
 {(a)!!(a-b+n-1)!!(b+n)!!}
} \quad ,
$$
where, $n!!:=1!2!3! \cdots n!$, and, of course, $n!:=1 \cdot 2 \cdots n$.
 
Indeed, let the left and right sides be $L_n (a,b)$ and $R_n (a,b)$
respectively. Dodgson's rule immediately implies that the recurrence:
$$
X_n(a,b)={{X_{n-1}(a,b) X_{n-1}(a+1,b+1)-X_{n-1} (a+1,b) X_{n-1} (a,b+1)}
\over {X_{n-2} (a+1,b+1)}} \quad,
$$
holds with $X=L$.
Since $L_n (a,b)=R_n (a,b)$ for $n=0,1$ (check!), and the recurrence
also holds with $X=R$ (check!\footnote{$^2$}
{\eightrm Divide both sides by the left, then use r!!/(r-1)!!=r!
whenever possible, and then r!/(r-1)!=r whenever possible, reducing
it to a completely routine polynomial identity.}
), 
it follows by induction that $L_n(a,b)=R_n(a,b)$ for {\it all} $n$. \halmos.
 
The special case $a=2n+1$, and $b=n$, reduces to a special case
of a conjecture of Enrico Bombieri\footnote{$^3$}
{\eightrm Who applies his imagination not only to mathematics, but also to art:
he is an accomplished painter.}, David Hunt, and
Alf van der Poorten ([BHP]). While this special case
was already done in [BHP] using a different method, I believe that
Dodgson's method should be extendable to prove their full conjecture.
\eject
{\bf References}
 
[BHP] E. Bombieri, D.C. Hunt, and A.J. van der Poorten, {\it Determinants
in the study of Thue's method and curves with prescribed singularities},
J. Experimental Mathematics, to appear.
 
[D] C.L. Dodgson, {\it Condensation of Determinants}, Proceedings
of the Royal Society of London {\bf 15}(1866), 150-155.
 
[M] P.A. MacMahon, {\it ``Combinatory Analysis''}, 
Cambridge University Press, 1918. [reprinted by Chelsea, 1984].
\bye